\documentclass[fleqn]{mat01}
\usepackage{times,mathtimy,amssymb,latexsym,amscd,epsfig}
\begin{document}

\makeatletter
\def\artpath#1{\def\@artpath{#1}}
\makeatother \artpath{C:/mathsci-arxiv/may2006/texfiles}

\setcounter{page}{147} \firstpage{147}

\def\d{\mbox{\rm d}}
\def\e{\mbox{\rm e}}

\newcommand{\latt}{\mbox{$\mathcal{L}$}}
\newcommand{\sltz}{\mbox{${\text{\rm SL}_2(\Z)}$}}
\newcommand{\suchthat}{\!:\;}
\newcommand{\Hy}{\mbox{${\mathbb H}$}}
\newcommand{\Z}{\mbox{${\mathbb Z}$}}
\newcommand{\R}{\mbox{${\mathbb R}$}}

\newtheorem{defi}{\rm DEFINITION}
\newtheorem{exam}{Example}
\newtheorem{propos}{\rm PROPOSITION}[ssection]
\newtheorem{corol}{\rm COROLLARY}[ssection]
\newtheorem{pop}{\it Proof of Proposition}[ssection]
\def\nota{\trivlist \item[\hskip \labelsep{{\it Notation}:}]}
\def\rem{\trivlist \item[\hskip \labelsep{\it Remark.}]}
\def\defini{\trivlist \item[\hskip \labelsep{\rm DEFINITION}]}

\title{Non-Euclidean visibility problems}

\markboth{Fernando Chamizo}{Non-Euclidean visibility problems}

\author{FERNANDO CHAMIZO}

\address{Departamento de Matem\'{a}ticas, Facultad de Ciencias,
Universidad Aut\'{o}noma \hbox{de Madrid,} 28049 Madrid, Spain}

\volume{116}

\mon{May}

\parts{2}

\pubyear{2006}

\Date{MS received 14 November 2005}

\begin{abstract}
We consider the analog of visibility problems in hyperbolic plane
(represented by Poincar\'{e} half-plane model $\Hy$), replacing the
standard lattice $\Z\times \Z$ by the orbit $z=i$ under the full modular
group $\sltz$. We prove a visibility criterion and study orchard problem
and the cardinality of visible points in large circles.\\

\end{abstract}

\keyword{Modular group; hyperbolic plane; Poincar\'{e} half-plane model.}

\maketitle

\section{Introduction}

Consider in $\R^2$ the standard lattice $L=\Z\times\Z$ and the origin
$(0,0)$. A point $(m,n)\in L-\{(0,0)\}$ is said to be {\it visible} if
the segment connecting the origin and $(m,n)$ does not contain any other
lattice points.

Visibility problems have been studied since a century. Perhaps the
most celebrated problems are the visible version of Gauss circle
problem and the so-called orchard problem (see other problems in
\cite{Erd}). In both of these problems one considers visible
lattice points in a large circle. The first problem consists of
approximating the cardinality of this set of points. It turns out
that improvements on the trivial bounds of the error term are
related to Riemann Hypothesis (see \cite{Now}). In the orchard
problem the visible points are considered to be thick and it is
asked the minimal thickness such that all exterior points are
eclipsed. In the formulation included in p.~150 of \cite{Pol},
``How thick must the trunks of the trees in a regularly spaced
circular orchard grow if they are to block completely the view
from the center?''. In contrast with the previous problem, orchard
problem can be considered as solved in a wide sense (see
\cite{All2}) by elementary methods.

In this paper we deal with the hyperbolic analog of visibility problems.
Namely, we consider Poincar\'{e}'s plane $\Hy$, i.e., the upper half
plane endowed with the metric
\begin{equation}
\d s^2=y^{-2}\ \d x^2+y^{-2}\ \d y^2, \label{metric}
\end{equation}
the origin $i\in \Hy$ and $\latt$ to be the orbit of $z=i$ under the
full modular group $\sltz$ (note that in the Euclidean case the lattice
$\Z\times \Z$ is the orbit of the origin under the discrete group formed
by all integral translations). We say that $z\in \latt$, $z\ne i$ is
{\it visible} if the arc of geodesic connecting $i$ and $z$ does not
contain any other point in $\latt$.

One cannot draw a parallel between the study of visibility
problems in a hyperbolic case and the Euclidean case due to the
following algebraic and geometric facts: Firstly, in the Euclidean
case the group of integral translations is Abelian, but in the
hyperbolic case the underlying group $\sltz$ is not. Secondly, the
Euclidean isoperimetric inequality $4\pi A\le l^2$, which is sharp
for circles, is qualitatively different from its hyperbolic analog
$4\pi A+A^2\le l^2$, for large areas (p.~11 of \cite{Iwa}).

We shall structure each of the following sections stating
Euclidean results first and then their hyperbolic counterparts;
this will ease the comparison between both settings. After
studying the symmetries and some other preliminary topics in \S2,
we give in \S3 a hyperbolic criterion for the visibility of a
point and investigate the structure of `lattice' points in rays.
Visible Gauss circle problem and orchard problem are discussed in
\S\S4 and 5. Finally, in \S6 we show some numerical data to
illustrate our results.

As an aside, we want to point out that although our main
motivation is number theoretical, visibility problems have called
the attention of some physicists (e.g., \cite{All1} and
\cite{Baa}) and it is plausible that the change of the geometry
could be meaningful in some applications. For instance, Olbers'
paradox (which even after two centuries still motivates some
research and controversy \cite{Wes}) in an idealized sharper form
taking into account occultation could lead to considerations about
visible points in some non-Euclidean space.

\begin{nota}
As usual, we shall represent with the same symbols a matrix in $\sltz$
and its associated fractional linear transformation acting on $\Hy$:
\begin{equation*}
\gamma=
\begin{pmatrix}
a&b\\ c&d
\end{pmatrix}
\longleftrightarrow
\gamma(z)=\frac{az+b}{cz+d}.
\end{equation*}
(Of course there is some ambiguity that would be avoided considering
$\textrm{PSL}_2(\Z)=\sltz /\{\pm I\}$.) With this convention we can
define the transpose of a fractional linear transformation or apply a
matrix to $z\in\Hy$.
\begin{equation*}
\gamma(z)=\frac{az+b}{cz+d}
\ \Rightarrow\
\gamma^t(z)=\frac{az+c}{bz+d},
\qquad
\begin{pmatrix}
a&b\\ c&d
\end{pmatrix}
z=\frac{az+b}{cz+d}.
\end{equation*}
We shall employ standard notation for the identity and symplectic
matrices, corresponding to identity function and involutive inversion
\begin{equation*}
I=\begin{pmatrix}
1&0\\ 0&1
\end{pmatrix}
\qquad\text{and}\qquad
j=
\begin{pmatrix}
0&-1\\ 1&0
\end{pmatrix}.
\end{equation*}
As we mentioned before, we shall let $\latt$ denote the orbit of $i$
under $\sltz$. We shall employ $\latt^*$ as an abbreviation of
$\latt-\{i\}$,
\begin{align*}
\latt &=\{\gamma(i)\!:\; \gamma\in\sltz \},\\[.3pc]
\latt^* &=\latt-\{i\}=\{\gamma(i)\! :\; \gamma\in\sltz ,\ \gamma\ne \pm I,\pm j \}.
\end{align*}
Finally, we shall denote by $d$ the hyperbolic distance corresponding to
Poincar\'e's metric (\ref{metric}).
\end{nota}

\section{Preliminaries and symmetries}

The group of proper motions of $\Hy$ is represented in $\text{\rm
SL}_2(\R)$ (in fact adding negative conjugation we get all motions), so
\begin{equation*}
d(\gamma(z),\tau(w))=d(z,\gamma^{-1}\tau(w)),\quad \forall z,w\in \text{\rm SL}_2(\R).
\end{equation*}
It is possible to write an explicit formula for the distance $d$ (see
\cite{Iwa} and \cite{Hux} for it and its geometrical interpretation)
that in the orbit of $z=i$ acquires an especially simple form (see
\cite{Cha})
\begin{align}
\gamma=
\begin{pmatrix}
a&b\\ c&d
\end{pmatrix}
\Rightarrow 2\; \text{cosh}\; d(i,\gamma(i))=a^2+b^2+c^2+d^2.
\label{dinl}
\end{align}
We find it convenient to state separately a calculation for further reference.

\begin{lemma}
Let $\gamma\in \text{\rm SL}_2(\R)$
\begin{align*}
\gamma=
\begin{pmatrix}
a&b\\ c&d
\end{pmatrix}
\Rightarrow
\gamma(i)=\frac{(ac+bd)+i}{c^2+d^2}
\qquad\text{and}\qquad |\gamma(i)|^2=\frac{a^2+b^2}{c^2+d^2},
\end{align*}
in particular{\rm ,} $(ac+bd)^2+1=(a^2+b^2)(c^2+d^2)$.
\end{lemma}

\begin{proof}
The formula for $\gamma(i)$ is just a straightforward computation, and
$(ac+bd)^2+1=(a^2+b^2)(c^2+d^2)$ can be quickly obtained from
$\text{det}(\gamma\gamma^t)=1$, noting that $a^2+b^2$ and $c^2+d^2$ are
the diagonal entries of $\gamma\gamma^t$ and the off-diagonal are
$ac+bd$.\hfill$\Box$
\end{proof}

It is clear that the set of Euclidean visible points has
eight-fold symmetry given by the dihedral group $D_4$. Namely {\it
in $\Z\times\Z${\rm ,} one of the eight points $(\pm x,\pm y)${\rm
,} $(\pm y,\pm x)${\rm ,} is visible if and only if the rest of
them are visible}. In the hyperbolic case, we have four-fold
symmetry.

\begin{lemma}
Let $z=x+iy\in \latt${\rm ,} then one of the points $z${\rm ,} $\overline{z}^{-1}${\rm ,}
$-{z}^{-1}${\rm ,} $-\overline{z}${\rm ,} is visible if and only if the rest of them
are visible.
\end{lemma}

\begin{proof}
The maps $T_k\!\!:\Hy\longrightarrow\Hy$ given by $T_1(z)=z$,
$T_2(z)=\overline{z}^{-1}$, $T_3(z)=-{z}^{-1}$, $T_4(z)=-\overline{z}$,
are involutive isometries with $T_k(i)=i$. They leave $\latt$ invariant
because
\begin{align*}
T_2\left(\begin{pmatrix} a&b\\ c&d\end{pmatrix}i\right)=
\begin{pmatrix} -c&d\\ -a&b\end{pmatrix}i,
\qquad
T_3(\gamma i)=
j\gamma(i)\quad\ \text{and}\ \quad
T_4=T_2\circ T_3.
\end{align*}
Hence if $g$ is an arc of geodesic with $i\in g$ and $\# (g\cap
\latt^*)=1$, then $T_k g$ has the same property.\hfill $\Box$
\end{proof}

This result allows to subdivide $\Hy$ in four `quadrants' with
disjoint interior:
\begin{align*}
Q_1 &=\{z\in \Hy\!: |z|\le 1,\  \text{Re}(z)\ge 0\},\quad
Q_2=\{z\in \Hy\!: |z|\ge 1,\  \text{Re}(z)\ge 0\},\\[.3pc]
Q_3&=\{z\in \Hy\!: |z|\ge 1,\  \text{Re}(z)\le 0\}, \quad
Q_4=\{z\in \Hy\!: |z|\le 1,\  \text{Re}(z)\le 0\}.
\end{align*}
With the notation of the previous proof
$T_k\big|_{Q_k}\!:Q_k\longrightarrow Q_1$ are well-defined isometries.

In the Euclidean case we can assign bijectively to each point in the
lattice the integral translation applying the origin on it, but in the
hyperbolic case the group $\sltz$ is not faithfully represented by the
orbit of $i$ due to the fact that the stability group of $z=i$ is $\{\pm
I, \pm j\}$. Hence the map
\begin{align*}
\sltz &\longrightarrow \latt=\{\gamma (i)\suchthat \gamma\in\sltz\}\\[.3pc]
\gamma &\longmapsto \gamma (i)
\end{align*}
is 4-to-1. We recover the Euclidean situation with some sign conventions.

\begin{lemma}
The map
\begin{equation*}
\left\{
\left(\begin{matrix}a&b \\ c&d\end{matrix}\right)\in\sltz\!:
\begin{array}{c}
a,b\ge 0,\ ac+bd> 0,\\ \ a^2+b^2<c^2+d^2
\end{array}
\right\}
\longrightarrow
\latt\cap \text{\rm Int}(Q_1)
\end{equation*}
given by $\gamma\mapsto \gamma (i)${\rm ,} is well-defined and bijective.
\end{lemma}

\begin{proof}
By Lemma~2.1, $ac+bd>0$ and $(a^2+b^2)/(c^2+d^2)<1$ is equivalent to
$\gamma(i)\in \text{\rm Int}(Q_1)$. On the other hand,
$\gamma(i)=\gamma'(i)$ if and only if $\gamma=\gamma'\tau$ with
$\tau\in\{\pm I,\pm j\}$ and there is only a choice of $\tau$ giving
$a,b\ge 0$.\hfill $\Box$
\end{proof}

\section{Rays and visible points}

Consider a half-infinite line $\ell$ starting at $(0,0)$ and
containing some other point of $L=\Z\times\Z$. It is fairly easy
to prove that {\it $\ell\cap L$ is the set of non-negative
multiples of a certain visible point $P\in\ell$}. On the other
hand, we have the straightforward arithmetic interpretation that
{\it visible points are simply the lattice points having coprime
coordinates}. In this section we shall state the hyperbolic analog
of these results.

\begin{propos}$\left.\right.$\vspace{.5pc}

\noindent Let $r$ be a ray in $\Hy$ {\rm (}half-infinite
geodesic{\rm )} starting at $i$ and containing some other point of
$\latt${\rm ,} then there exists $\gamma\in\sltz$ such that
\begin{equation*}
r\cap \latt^*=
\big\{(\gamma\gamma^t)^{n-1}\gamma(i)\suchthat n\in \Z^+\big\}\cup
\big\{(\gamma\gamma^t)^{n}(i)\suchthat n\in \Z^+\big\}.
\end{equation*}
Moreover $\gamma(i)$ is visible and the points in $r\cap\latt$ are equally spaced on $r$.
\end{propos}

The proof employs the following auxiliary result:

\setcounter{defin}{1}
\begin{lemma}
Let $\gamma\in\sltz$ with $\gamma(i)\ne i${\rm ,} then $\tau=\gamma\gamma^t$
leaves invariant the geodesic connecting $i$ and $\gamma(i)${\rm ,} and it
holds $d(i,\gamma(i))=d(\gamma(i),\tau(i))$.
\end{lemma}

\begin{proof}
A (non-vertical infinite) geodesic $g$ can be represented as a
Euclidean semicircle in $\Hy$ orthogonal to the real axis. If
$i\in g$, by simple trigonometry, end-points are $-x_0$ and
$x_0^{-1}$ for some $x_0\in\R$.

Take the geodesic passing through $i$ and $\gamma(i)$ as $g$, then
$i\in \gamma^{-1}g$. Hence we have $\gamma^{-1}(-x_0)=-y_0$ and
$\gamma^{-1}(x_0^{-1})=y_0^{-1}$ for some $y_0\in\R$. Consequently
the end-points of $\gamma j^{-1}\gamma^{-1}j g$ are $\gamma
j^{-1}\gamma^{-1}j(-x_0)=-x_0$ and $\gamma
j^{-1}\gamma^{-1}j(x_0^{-1})=x_0^{-1}$. It means that $\gamma
j^{-1}\gamma^{-1}j$ lea\-ves $g$ invariant. A calculation proves
$j^{-1}\gamma^{-1}j=\gamma^t$.

On the other hand, $d(i,\gamma^{t}(i))=d(\gamma(i),\tau(i))$ because
$\gamma$ is an isometry, and $d(i,\gamma(i))=d(i,\gamma^t(i))$ by
(\ref{dinl}).\hfill $\Box$
\end{proof}

\begin{pop}
{\rm Parametrizing $r$ by arc length, each point on $r$ is determined by its
distance to $i$. Plainly $r$ contains exactly one visible point, say
$\gamma(i)$ and write $l=d(i,\gamma(i))$.

Note that the signs of $\text{Re}(\gamma(i))$ and
$\text{Re}(\gamma\gamma^t(i))$ coincide (see Lemma~2.1). Then by the
previous lemma $\gamma\gamma^t$ applies the half geodesic $r$ into
itself and
\begin{align*}
d(i,\gamma\gamma^t(i))&= d(i,\gamma(i))+d(\gamma(i),\tau(i))\\[.3pc]
&=2d(i,\gamma(i))=2l.
\end{align*}
Hence $\gamma\gamma^t\big|_r$ is just a translation of lenght $2l$ along $r$ and
\begin{align*}
\big\{(\gamma\gamma^t)^{n-1}\gamma(i)\suchthat n\in \Z^+\big\}&\cup
\big\{(\gamma\gamma^t)^{n}(i)\suchthat n\in \Z^+\big\}\\[.2pc]
&=\{z\in r\!:\; d(i,z)\in l\Z^+\}.
\end{align*}
This set of $l$-spaced points is contained in $r\cap \latt^*$. It
remains to prove that any $w\in r\cap \latt^*$ satisfies
$d(i,w)\in l\Z^+$. If $d(i,w)\not\in l\Z^+$ then for some
$k\in\Z$,
\begin{equation*}
2kl<d(i,w)<(2k+1)l
\quad\text{or}\quad
(2k-1)l<d(i,w)<2kl.
\end{equation*}
In the first case $z=(\gamma\gamma^t)^{-k}(w)\in r\cap \latt^*$ verifies
$0\ne d(i,z)<l=d(i,\gamma(i))$ which is a contradiction because
$\gamma(i)$ is visible. In the second case the same argument applies
with $z=j(\gamma\gamma^t)^{-k}(w)$. In this connection note that
$(\gamma\gamma^t)^{-k}(w)\in r'$ where $r'$ is the complementary ray of
$r$, i.e. $r\cap r'=\{i\}$ and $r\cup r'$ form an infinite geodesic $g$;
and $j$ applies $r'$ isometrically into $r$ leaving $i$ fixed ($j$
permutes the end-points of $g$), then $d(i,z)=d(i,
(\gamma\gamma^t)^{-k}(w))=2kl-d(i,w)$.\hfill $\Box$}
\end{pop}

Now we are going to characterize visible points in terms of their
coordinates. Recall firstly that any $z\in\latt^*$ is uniquely written
as
\begin{equation*}
z=\frac{B+i}D\qquad\text{with}\,\,B,D\in\Z,
\end{equation*}
and consider the map
\begin{equation*}
\begin{array}{rcl}
&V\!:\latt^* \longrightarrow  L^*=\Z\times\Z-\{(0,0)\}\\[.3pc]
&z=(B+i)/D \longrightarrow  (B,D-A)
\end{array}
\qquad\text{with}\quad A=(B^2+1)/D.
\end{equation*}
Note that it is well-defined and applies the first quadrant $\latt^*\cap
Q_1$ into the Euclidean first quadrant (see \S2, esp. Lemma~2.1). By
symmetry, it is enough to state visibility criterion in $Q_1$; in the
rest of the quadrants it is similar up to sign changes.

\begin{theorem}[\!]
Let $z\in\latt^*\cap Q_1$ and $A,B,D$ as before{\rm ,} then $z$ is not visible
if and only if there exists integers $1\le a\le b<d$ with $ad=b^2+1$ and
$b|B${\rm ,} such that
\begin{equation*}
\frac Bb=\frac{D-A}{d-a}\ne 1.
\end{equation*}
\end{theorem}

\begin{proof}
Firstly note that (see Lemma~2.1)
\begin{equation*}
z=\gamma(i)\ \Rightarrow\ \gamma\gamma^t=
\begin{pmatrix}
A&B\\ B&D
\end{pmatrix}.
\end{equation*}
By Lemma~3.2, the hyperbolic motion $\gamma\gamma^t$ leaves invariant
the geodesic $g$ connecting $i$ and $\gamma(i)$. End-points of $g$, say
$z_1,z_2\in\R$, are the roots of the quadratic equation
$(Az+B)/(Bz+D)=z$, then the Euclidean center of the semicircle
representing $g$ is
\begin{equation}
\frac{z_1+z_2}2=\frac{A-D}{2B}.
\label{center}
\end{equation}
Now we shall consider both implications separately.

$(\Rightarrow )$ If $z=\gamma(i)$ is not visible, let $\tau(i)\ne
z$ be the visible point in the ray connecting $i$ and $z$. Let us
take $\gamma$ and $\tau$ normalised as in Lemma~2.3. By
eq.~(\ref{center}) and Lemma~3.2
\begin{equation*}
\tau\tau^t =
\begin{pmatrix}
a&b\\ b&d
\end{pmatrix}
\ \Rightarrow\ \frac{a-d}{2b}=\frac{A-D}{2B}
\end{equation*}
because $\tau(i)$ and $\gamma(i)$ are on the same ray. Of course
$\tau\tau^t\in\sltz$ $\Rightarrow$ $ad=b^2+1$ and it only remains to prove
$b|B$, $b\ne B$. By Proposition~3.1, $\gamma=(\tau\tau^t)^n\tau$ or
$\gamma=(\tau\tau^t)^n$ with $n\in\Z^+$. In any case
$\gamma\gamma^t=(\tau\tau^t)^k$, $k\ge 2$ and by induction on $k$ it
follows that $b$ divides the second entry of $(\tau\tau^t)^k$. The
positivity given by Lemma~2.3 assures $b<B$.

$(\Leftarrow )$ As $ad=b^2+1$, by the theory of binary quadratic
forms (see Art.\;183, \cite{Gau}), we can find $\tau\in\sltz$ such
that
\begin{equation*}
\tau\tau^t=
\begin{pmatrix}
a&b\\ b&d
\end{pmatrix}.
\end{equation*}
In fact we can assume that $\tau$ is as in Lemma~2.3. The relation
$B/b=(D-A)/(d-a)$ and (\ref{center}) imply that $\tau(i)$ and
$\gamma(i)$ are in the same geodesic ray starting at $i$. Since $B>b$
and $D-A>d-a>0$,
\begin{align*}
\left.
\begin{matrix}
D^2+A^2-2AD>d^2+a^2-2ad
\\[.3pc]
4(B^2+1)>4(b^2+1)
\end{matrix}
\right\}
&\Rightarrow\
(D+A)^2>(d+a)^2\\[.2pc]
&\Rightarrow\
D+A>d+a
\end{align*}
just adding both inequalities. Hence
\begin{align*}
D-A>d-a,\quad D+A>d+a\ \Rightarrow\ D>d \ \Rightarrow\ \text{Im}(\gamma(i))<\text{Im}(\tau(i)).
\end{align*}
As $\gamma(i)$ and $\tau(i)$ belong to the same ray, the latter
condition implies
\begin{equation*}
d(i,\gamma(i))>d(i,\tau(i)),
\end{equation*}
thus $\gamma(i)$ is not visible.\hfill$\Box$
\end{proof}

\begin{exam}
{\rm Consider
\begin{equation*}
z=\begin{pmatrix} 2&1\\ 3&2\end{pmatrix}i=\frac{8+i}{13}
\ \Rightarrow\ B=8,\ D=13,\ A=5,\ D-A=8.
\end{equation*}
In this case, the conditions of the theorem read $b=1,4$, $ad=5,17$,
respectively, with $b=a-d\ne 8$. This is fulfilled for $b=1$, $d=5$,
$a=1$, and hence $z=(8+i)/13$ is not visible.}
\end{exam}

\begin{exam}
{\rm The point
\begin{equation*}
z=\begin{pmatrix} 1&3\\ 2&7\end{pmatrix}i=\frac{23+i}{53}
\end{equation*}
is visible because $B=23$, $D=53$, $A=10$ and as $B$ and $D-A$ are
coprime the equation $B/b=(D-A)/(d-a)\ne 1$ cannot hold.}
\end{exam}\pagebreak

Let us state separately the last remark:

\setcounter{corol}{3}
\begin{corol}$\left.\right.$\vspace{.5pc}

\noindent If $B$ and $D-A$ are coprime then $z$ is visible. Equivalently{\rm ,} if
$V(z)$ is visible {\rm (}in Euclidean sense{\rm )} then $z$ is visible {\rm (}in
hyperbolic sense{\rm )}.
\end{corol}

\begin{rem}
The reciprocal is not true. The simplest counterexample is $z=(2+i)/5$
which is visible with $V(z)=(B,D-A)=(2,4)$.
\end{rem}

In the Euclidean case, if we enumerate the points on each ray starting
by zero (assigned to the origin), then the points labelled by even
numbers, say the {\it points in even place}, are exactly the sublattice
of points with even coordinates. In the hyperbolic case we can define
equally points in even place and the following result allows a
coordinate characterization.

\setcounter{propos}{4}
\begin{propos}$\left.\right.$\vspace{.5pc}

\noindent $z\in\latt$ is a point in even place if and only if
$z=\tau(i)$ with $\tau$ symmetric and equivalently{\rm ,} if and only if
\begin{equation*}
z=\frac{(a+d)b+i}{b^2+d^2}
\end{equation*}
for some integers $ad=b^2+1$.
\end{propos}

\begin{proof}
By Proposition~3.1 we have that the points in $r\cap\latt$, where $r$ is
the ray connecting $i$ and $z$, are equally spaced. In fact we have
proved that $\gamma\gamma^t\big|_r$ is a translation of length $2l$
where $l$ is the separation between consecutive points. Hence $z$ is in
even place if and only if $\tau=(\gamma\gamma^t)^n$ where $n$ is a
non-negative integer. Then if $z$ is in even place, $\tau$ is symmetric.
Reciprocally, if $\tau$ is symmetric we can write (Art.\;183,
\cite{Gau}) (as in the previous proof) $\tau=\delta\delta^t$, and by
Proposition~3.1, $\delta=(\gamma\gamma^t)^{n-1}\gamma$ or
$\delta=(\gamma\gamma^t)^{n}$. In any case, $\tau=(\gamma\gamma^t)^m$
and $z$ is in even place.\hfill $\Box$
\end{proof}

\section{The visible lattice point problem}

The asymptotics of the number of visible points in a Euclidean circle
of large radius $R$ has been studied by several authors. This number is
usually approximated by a formula like
\begin{equation}
E^*(R)=\frac{6}{\pi^2}R^2+O(R^\alpha).
\label{eq:4eu}
\end{equation}
Surprisingly, any improvement on the trivial exponent $\alpha=1$ (see
\cite{Now}) lead to considerations on Riemann Hypothesis and there are
no unconditional results with $\alpha<1$. Several authors have proved
(\ref{eq:4eu}) for some $\alpha$ assuming Riemann Hypothesis (using the
arguments in \cite{Now} and intricate exponential sums estimates). The
best conditional result so far is $\alpha=221/608+\epsilon$ for every
$\epsilon>0$ \cite{Wu}.

In hyperbolic setting, the relation with Riemann Hypothesis
disappears, rough\-ly speaking because most of the lattice points
stay close to the boundary and hence comparatively few points are
eclipsed, and the contribution of invisible points is absorbed by
error term. Considering firstly all the points in the orbit of
$i$, visible and invisible, the asymptotics of the number of
points in a large circle of radius $R$ is given as (see
\cite{Phi}, we introduce a $1/4$ extra factor because $4$ is the
cardinality of stability group of $i$)\pagebreak
\begin{equation}
H(R)=\frac 32 \e^R+O(\e^{\alpha R}).
\label{eq:4hy}
\end{equation}
Using harmonic analysis on $\Hy\big\backslash\sltz$ one can get
$\alpha=2/3$ (see \cite{Iwa}, \S12). This exponent has not been
improved but the natural conjecture (supported by average results
\cite{Cha}) is $\alpha=1/2+\epsilon$ for every $\epsilon>0$.

\setcounter{propos}{0}
\begin{propos}$\left.\right.$\vspace{.5pc}

\noindent Let $H^*(R)$ be the number of visible points in the circle
$\{z\in\Hy\!:\; d(i,z)\le R\}${\rm ,} then
\begin{equation*}
H^*(R)=H(R)-\frac{3}{2} \e^{R/2}+O(\e^{R/3}).
\end{equation*}
\end{propos}

\begin{proof}
Given a  ray $r$ connecting $i$ and some point in $\latt^*$, let
\begin{equation*}
r(R)=\#\{z\in\latt^*\cap r\!:\; d(i,z)\le R\}.
\end{equation*}
According to Proposition~3.1, the points in $\latt\cap r$ are
$l$-spaced, hence
\begin{equation*}
r(R)=\left[\frac Rl\right]\qquad
\text{and}\qquad
1=\sum_{n\le R}\mu(n)r(R/n)\quad\text{ for }R\ge l,
\end{equation*}
(see \cite{Ell}) where $[\;\cdot\; ]$ denotes integral part and $\mu$ is
M\"obius function.

Let $\mathcal{R}$ be the set of rays as before containing some
$z\in\latt^*$ with $d(i,z)\le R$. Each ray in $\mathcal{R}$
contains exactly a visible point and we have for $R>1$,
\begin{align*}
H^*(R)=\sum_{r\in\mathcal{R}}1&=\sum_{n\le R}\mu(n)\sum_{r\in\mathcal{R}}r(R/n)\\[.6pc]
&=\sum_{n\le R}\mu(n)H(R/n).
\end{align*}
Using (\ref{eq:4hy})
\begin{equation*}
H^*(R)=H(R)-H(R/2)-H(R/3)+O(\e^{\alpha R/5})
\end{equation*}
and taking $\alpha=2/3$, $H(R/2)-H(R/3)=3\e^{R/2}/2+O(\e^{R/3})$
we get the result. Note that under the conjecture
$\alpha=1/2+\epsilon$ we could diminish error term to
$O(\e^{(1+\epsilon)R/4})$ extracting an extra $-3\e^{R/3}/2$
term.\hfill$\Box$
\end{proof}

The last proposition allows to translate to $H^*(R)$ some  results known for $H(R)$.
Following \cite{Phi}, we define the normalized remainder
\begin{equation*}
\Delta^*(R)=\frac{H^*(R)-3\e^{R}/2}{\e^{R/2}}.
\end{equation*}
It turns out that $\Delta^*(R)$ is biased (because of the
influence on invisible points), and inherits the oscillation of
$H(R)$. After Proposition~4.1, this is just a consequence of the
main results of \cite{Phi}. \pagebreak

\setcounter{corol}{1}
\begin{corol}$\left.\right.$\vspace{.5pc}

\noindent The mean value of $\Delta^*(R)$ is $3/2${\rm ,} i.e.
\begin{equation*}
\lim_{R\to\infty}\frac 1R\int_1^R\Delta^*(t)\; \d t=3/2,
\end{equation*}
but $\Delta^*(R)$ is not bounded. In fact
\begin{equation*}
\limsup_{R\to \infty}
\frac{\Delta^*(R)}{(\log R)^{\delta}}=\infty\qquad
\text{for every }\quad \delta<1/4.
\end{equation*}
\end{corol}

\begin{proof}
Let $\Delta(R)=(H(R)-3\e^R/2)/\e^{R/2}$. By Theorems~1.1 and 1.2
in \cite{Phi} (note that for the full modular group
$E(z,s)=\zeta_Q(s)/\zeta(2s)$ where $\zeta_Q$ is an Epstein zeta
function \cite{Hux}), it holds that
\begin{equation*}
\lim_{R\to\infty}\frac 1R\int_1^R\Delta(t)\; \d t=0
\quad\text{and}\quad \limsup_{R\to \infty} \frac{\Delta(R)}{(\log
R)^{\delta}}=\infty.
\end{equation*}
By Proposition~4.1, $\Delta^*(R)=\Delta(R)-3/2+O(\e^{-R/6})$ and
the claimed results follow.

\hfill$\Box$\vspace{-.5pc}
\end{proof}

We have also some control on a quantity related to the variance.

\begin{corol}
\begin{equation*}
\limsup_{R\to \infty}
\frac 1R\int_1^R|t^{-1}\Delta^*(t)|^2\; \d t<\infty.
\end{equation*}
\end{corol}

\begin{proof}
Corollary~2.1.1 of  \cite{Cha} implies
\begin{equation*}
\int_{X}^{2X}\left\vert H\left(\text{\rm arc\;cosh}\; \frac{x}{2}\right)-\frac{3}{2} x\right\vert^2\; \d x=O(X^2\log^2 X).
\end{equation*}
With the change of variables $x=2\; \text{\rm cosh}\; t$ and
writing $r=\log X$, we have
\begin{equation*}
\int_{r}^{r+1}\bigg|H(t)-\frac 32 \e^t\bigg|^2\e^t\; \d
t=O(r^2\e^{2r}).
\end{equation*}
Using Proposition~4.1, and after some manipulations we get that
\begin{equation*}
\int_r^{r+1}\left\vert \frac{\Delta^*(t)}{r}\right\vert ^2\; \d t
\end{equation*}
is bounded. Summing on $1\le r\le R-1$, the result is proved.\hfill$\Box$
\end{proof}

\section{The orchard problem}

It is known that the solution of the orchard problem (as stated in the
Introduction) is that the view is obstructed in a circular orchard of
radius $R$ if and only if the trunks have radii $\epsilon\ge
R^{-1}+f(R)$ for certain $f(R)=O(R^{-2})$ (an \lq exact' formula is
given in \cite{All2}). One can also consider, so to speak, negative
orchard problem, asking for the maximal $\epsilon$ such that it is
possible to see all the visible points in the circle of radius $R$.
In Euclidean setting the solution is the same as that of the original
problem, but in hyperbolic setting both problems considerably differ.

We shall associate to each $z\in\latt^*$ and $\epsilon>0$ the thick
point $z_{\epsilon}$ with radius $\epsilon$; this means the circle
$z_\epsilon=\{w\in\Hy\!:\; d(w,z)\le\epsilon\}$.

\begin{defini}$\left.\right.$\vspace{.5pc}

\noindent Given $z,w\in\latt^*$ we say that $z_\epsilon$ {\it eclipses} $w$ if
$r\cap z_{\epsilon}\ne\emptyset$ where $r$ is the geodesic arc
connecting $i$ and $w$.
\end{defini}

Our main tool for treating the obstruction of view in $\latt$ is
the following result.

\setcounter{propos}{0}
\begin{propos}$\left.\right.$\vspace{.5pc}

\noindent Let $z,w\in\latt^*\cap Q_1$ with $d(i,z)\le d(i,w)${\rm ,} say
$z=\gamma(i)$ and $w=\tau(i)$. Then $z_\epsilon$ eclipses $w$ if and
only if
\begin{equation*}
\text{\rm sinh}\;\epsilon\ge
\frac{\big|\text{\rm Tr}(\gamma\gamma^tj\tau\tau^t)\big|}{2\;\text{\rm sinh}\;d(i,w)},
\end{equation*}
where $\text{\rm Tr}$ indicates the  trace.
\end{propos}

\begin{proof}
Note firstly that $z_\epsilon$ eclipses $w$ if and only if $\epsilon$ is
greater than the distance of $z$ to the geodesic $g$ connecting $i$ and
$w$, because the foot of the perpendicular through $z$, say $F$, belongs
to $Q_1$ and $d(i,F)\le d(i,z)\le d(i,w)$ (by hyperbolic Pythagorean
theorem $\text{\rm cosh}\;a\; \text{\rm cosh}\; b=\text{\rm cosh}\;c$).

Whence we are going to prove that for every
$\gamma,\tau\in\text{SL}_2(\R)$, $z=\gamma(i)$, $w=\tau(i)\ne i$, if $g$
is the (infinite) geodesic through $i$ and $w$, we have
\begin{equation}
\text{\rm sinh}\; d(z,g)=
\frac{\big|\text{\rm Tr}(\gamma\gamma^tj\tau\tau^t)\big|}{2\;\text{\rm sinh}\;d(i,w)}.
\label{distpg}
\end{equation}
Consider $m_\theta\in\text{SL}_2(\R)$ given by
\begin{equation*}
m_\theta=\begin{pmatrix}\cos\theta&-\sin\theta\\ \sin\theta &\cos\theta \end{pmatrix}.
\end{equation*}
In hyperbolic plane this is a rotation at $i$ of angle $2\theta$ (\S1.2,
\cite{Iwa}). Both sides in (\ref{distpg}) are invariant under the
changes $\gamma\mapsto m_\theta\gamma m_{\theta'}$ and $\tau\mapsto
m_\theta\tau m_{\theta'}$ because $m_\theta,m_{\theta'}$ are orthogonal
matrices, leave $i$ invariant and $m_{\theta}^tjm_{\theta}=j$. With a
suitable choice of $m_\theta$ and $m_{\theta'}$ we can assume by
Cartan's decomposition (\S1.3, \cite{Iwa}) that
\begin{equation*}
\tau=
\begin{pmatrix}
\lambda &0
\\
0 &\lambda^{-1}
\end{pmatrix}
\quad \text{and}\quad z=\lambda^2 i\quad\text{for some }\ \lambda\in\R^+.
\end{equation*}
Let $\gamma(u)=(au+b)/(cu+d)$, then $\big|\text{\rm
Tr}(\gamma\gamma^tj\tau\tau^t)\big|=|ac+bd|(\lambda^2-\lambda^{-2})$ and
using (\ref{dinl}), we have that (\ref{distpg}) reduces to prove that
the hyperbolic distance $D$ from $z=\gamma(i)$ to the imaginary axis is
given by $\text{\rm sinh}\;D=|ac+bd|$. As $z$ is in the circle
$|\zeta|=|\gamma(i)|$ which is orthogonal to this axis, by hyperbolic
Pythagorean theorem
\begin{equation*}
\text{\rm sinh}^2\!\,D =\frac{\text{\rm
cosh}^2\;d(i,\gamma(i))}{\text{\rm cosh}^2\;d(i,|\gamma(i)|i)}-1,
\end{equation*}
which using (\ref{dinl}), Lemma~2.1 and   (\ref{metric}), gives the result.\hfill$\Box$
\end{proof}

First let us consider the negative orchard problem.

\begin{propos}$\left.\right.$\vspace{.5pc}

\noindent Let $C_R^*=\latt^*\cap\{z\!:\;d(i,z)\le R\}$. If $\epsilon<2\e^{-R}$ then
none of the points in $C_R^*$ enlarged to radius $\epsilon$ eclipses
another point in $C_R^*$.
\end{propos}

\begin{proof}
Let $z,w\in C_R^*\cap Q_1$, say $d(i,z)\le d(i,w)$, and let $d_1$ and $d_2$ be the distances of
$z$ and $w$ to the geodesics connecting $i$ with $w$ and $i$ with $z$, respectively. Using sine rule \cite{Hux}
\begin{equation*}
\frac{\text{\rm sinh}\;d_1}{\text{\rm sinh}\;d(i,z)}
=\frac{\text{\rm sinh}\;d_2}{\text{\rm sinh}\;d(i,w)},
\end{equation*}
hence $d_1\le d_2$ and $w_\epsilon$ eclipses $z$ implies that $z_\epsilon$ eclipses $w$.

It is easy to check that for a symmetric matrix in $\sltz$ the
off-diagonal entry and the trace are congruent modulo $2$. If $B$ and
$\beta$ are the off-diagonal entries of $\gamma\gamma^t$ and
$\tau\tau^t$, a calculation proves
\begin{align*}
\text{Tr}(\gamma\gamma^tj\tau\tau^t) &\equiv  B\text{Tr}(\tau\tau^t)+\beta\text{Tr}(\gamma\gamma^t)\qquad (\text{mod }2),
\\[.3pc]
&\equiv 2\text{Tr}(\tau\tau^t)\text{Tr}(\gamma\gamma^t)\equiv 0\qquad (\text{mod }2).
\end{align*}
By Proposition~5.1, if $z_\epsilon$ eclipses $w$ then $\text{\rm sinh}\,
\epsilon\ge (\text{\rm sinh}\, R)^{-1}$, and this implies $\epsilon\ge
2\e^{-R}$.

If against our assumption $z$ and $w$ do not belong to $Q_1$ but the
corresponding rays determine an acute angle, then the same proof applies
after a suitable rotation. If the angle is not acute, if $z_\epsilon$
eclipses to $w$ then $i\in z_{\epsilon}$ and, according to (\ref{dinl}),
$\text{\rm cosh}\,\epsilon \ge 3/2$, i.e. $\epsilon\ge 0.9624\dots$ and
we can assume $2\e^{-R}<0.764\dots$ because otherwise $C^*_R=\emptyset$.\hfill $\Box$
\end{proof}

A construction using the properties of Fibonacci numbers allows to show
that the previous result is sharp.

\begin{propos}$\left.\right.$\vspace{.5pc}

\noindent Given $C>2$ there exist sequences of values
$z\in\Hy${\rm ,} $w\in\Hy$ and $R\in\R$ tending to $\infty$ such
that $z,w\in C_R^*$ are visible points and $z_\epsilon$ eclipses
$w$ with $\epsilon=C\e^{-R}$.
\end{propos}

\begin{proof}
Consider Fibonacci sequence
$\{F_n\}_{n=1}^\infty=\{1,1,2,3,5,8,\dots\}$ and for each $n$,
\begin{equation*}
\gamma=
\begin{pmatrix}
F_{6n-1} & F_{6n-2}
\\
F_{6n+1} & F_{6n}
\end{pmatrix}
\quad\text{and}\quad
\tau=
\begin{pmatrix}
F_{6n-1} & F_{6n}
\\
F_{6n+1} & F_{6n+2}
\end{pmatrix}.
\end{equation*}
It holds that $\gamma,\tau\in\sltz$ (use that $F_{k+1}/F_k$ are the
convergents of the golden ratio, or employ the recurrence formula,
4.2.3(d) and 4.3.9(b) of \cite{Lov}). Choose $z=\gamma(i)$, $w=\tau(i)$
and $R=d(i,w)$. By Lemma~2.3, $z,w\in\latt^*\cap Q_1$ and we are under
the hypothesis of Proposition~5.1.

Using the properties of Fibonacci numbers (see (4.2) and 4.2.3(d) of \cite{Lov}) we get
\begin{equation*}
\gamma\gamma^t=
\begin{pmatrix}
F_{12n-3} & F_{12n-1}
\\
F_{12n-1} & F_{12n+1}
\end{pmatrix}
\quad\text{and}\quad
\tau\tau^t=
\begin{pmatrix}
F_{12n-1} & F_{12n+1}
\\
F_{12n+1} & F_{12n+3}
\end{pmatrix}.
\end{equation*}
Take $m=12n-3$ or $12n-1$. By Euclidean algorithm $F_{m+1}$ and $F_m$ are
coprime, and $F_{m+1}+F_m$ and $F_{m+1}-F_m$ are coprime too (both are
odd numbers). Then $F_{m+2}=F_{m+1}+F_m$ and
$F_{m+4}-F_m=2F_{m+2}+F_{m+1}-F_m$ are also coprime. By Corollary~3.4,
$z$ and $w$ are visible.

A calculation shows
\begin{align*}
\text{Tr}(\gamma\gamma^t j\tau\tau^t)
&=
F_{12n-1}(F_{12n-1}-F_{12n+3})
+F_{12n+1}(F_{12n+1}-F_{12n-3})
\\[.2pc]
&=
-F_{12n-1}(F_{12n}+F_{12n+2})
+F_{12n+1}(F_{12n-2}+F_{12n})
\\[.2pc]
&=
(F_{12n+1}F_{12n-2}-F_{12n-1}F_{12n})\\
&\quad\,-(F_{12n+2}F_{12n-1}-F_{12n}F_{12n+1})
\\[.2pc]
&=
-1-1=-2.
\end{align*}
Where we have firstly used that $F_{k+2}-F_{k-2}=F_{k+1}+F_{k-1}$ and
secondly, as before, that $F_{k+1}/F_k$ are the convergents of
$(1+\sqrt{5})/2$.

By Proposition~5.1 we have that for $\text{\rm sinh}\; \epsilon\ge
(\text{\rm sinh}\; R)^{-1}$, $z_\epsilon$ eclipses $w$, and this
inequality holds with $\epsilon=C\e^{-R}$ for large enough $R$.\hfill$\Box$
\end{proof}

It turns out (see the proof below) that the unique way of blocking
completely the view from the origin in a hyperbolic orchard is
enlarging a certain fixed quantity the trunks of the first four
trees. In this sense, orchard problem becomes trivial in its
original form.

\begin{propos}$\left.\right.$\vspace{.5pc}

\noindent Every $w$ with $d(i,w)>R$ is eclipsed by some $z_\epsilon$
with $z\in C_R^*$ if and only if $\epsilon\ge \log(1+\sqrt{2})$.
\end{propos}

\begin{proof}
If $z=\gamma(i)=(ai+b)/(ci+d)$, $z\ne i$, we have shown at the end of
the proof of Proposition~5.1 that the distance $D$ from $z$ to the
imaginary axis verifies $\text{\rm sinh}\,D=|ac+bd|$, hence $\text{\rm
sinh}\; D\ge 1$ and we cannot block the view along the imaginary axis if
$\epsilon<\text{\rm arc\,sinh}\; 1=\log(1+\sqrt{2})$.

On the other hand, let $z_2=1+i\in Q_2\cap\latt^*$. The circle
$\{w\! :\; d(z_2,w)\le \log(1+\sqrt{2})\}$ correspond to the
Euclidean circle in $\Hy$ given by $(x-1)^2+(y-\sqrt{2})^2\le 1$
(see \S1.1, \cite{Iwa}). Applying $T_k^{-1}T_2$ we get four
intersecting circles around the origin blocking the view from
$z=i$.\hfill$\Box$
\end{proof}

Even disregarding near points, if we argue heuristically thinking
that the points in $ C_R^*$ are uniformly distributed along the
boundary (of length $2\pi\,\text{\rm sinh}\,R$), we can expect
maximal spacing as large as $\e^{R/2}$ (in particular unbounded,
in contrast with the Euclidean case). This effectively happens
when we pass from a quadrant to another. For instance, the rays
$r_{-},\;r_{+}$ connecting $i$ and $i-n,\; i+n$ are consecutive in
the circle $C_R^*$ where $\text{\rm cosh}\,R=(n^2+2)/2$ and the
spacing $d(i-n,r_{+})=d(i+n,r_{-})$ is comparable to $2\e^{R/2}$.
Applying elements of $\sltz$ the same phenomenon repeats at
different scales inside each quadrant.

\section{Numerical results}

\begin{table}[b]
\processtable{}
{\begin{tabular}{@{}c@{\hskip 1cm}c@{\hskip 1cm}c@{\hskip 1cm}c@{}}\hline
  & & &\\[-.7pc]
$\e^R$ & Visible   & Invisible  & Error\\\hline
  & & &\\[-.7pc]
\phantom{0}1000 & \phantom{0}1436 & \phantom{0}60 & $-16.56$ \\[.2pc]
\phantom{0}2000 & \phantom{0}2904 & \phantom{0}92 & $-28.91$ \\[.2pc]
\phantom{0}3000 & \phantom{0}4408 & 100 & \phantom{0}$-$9.84 \\[.2pc]
\phantom{0}4000 & \phantom{0}5960 & 124 & \phantom{$-$}54.86 \\[.2pc]
\phantom{0}5000 & \phantom{0}7336 & 140 & $-57.93$ \\[.2pc]
\phantom{0}6000 & \phantom{0}8844 & 148 & $-39.81$ \\[.2pc]
\phantom{0}7000 & 10372 & 160 & \phantom{0}$-$2.50 \\[.2pc]
\phantom{0}8000 & 11792 & 176 & $-73.83$ \\[.2pc]
\phantom{0}9000 & 13280 & 176 & $-77.69$ \\[.2pc]
10000 & 14880 & 184 &\phantom{$-$}30.00\\\hline
\end{tabular}}{}
\end{table}

Using Theorem~3.3 it is easy to write a computer program distinguishing
visible from invisible points in a large hyperbolic circle. The
cumulative number of them is given in table~1

\begin{table}[t]
\processtable{}
 {\begin{tabular}{@{}c@{\hskip 1cm}c@{\hskip 1cm}c@{}}\hline
 & &\\[-.7pc]
$\e^R$ & Invisible  & {Approx.} \\\hline
 & &\\[-.7pc]
\phantom{0}1000 & \phantom{0}60 & \phantom{0}63.66\\[.2pc]
\phantom{0}2000 & \phantom{0}92 & \phantom{0}87.52 \\[.2pc]
\phantom{0}3000 & 100 & $105.53$ \\[.2pc]
\phantom{0}4000 & 124 & $120.58$ \\[.2pc]
\phantom{0}5000 & 140 & $133.75$ \\[.2pc]
\phantom{0}6000  & 148 & $145.60$ \\[.2pc]
\phantom{0}7000  & 160 & $156.45$ \\[.2pc]
\phantom{0}8000  & 176 & $166.51$ \\[.2pc]
\phantom{0}9000  & 176 & $175.93$ \\[.2pc]
10000 & 184 & $184.82$ \\\hline
\end{tabular}}{}\vspace{-.7pc}
\end{table}

\begin{figure}[t]
\centerline{\epsfxsize=12.3cm\epsfbox{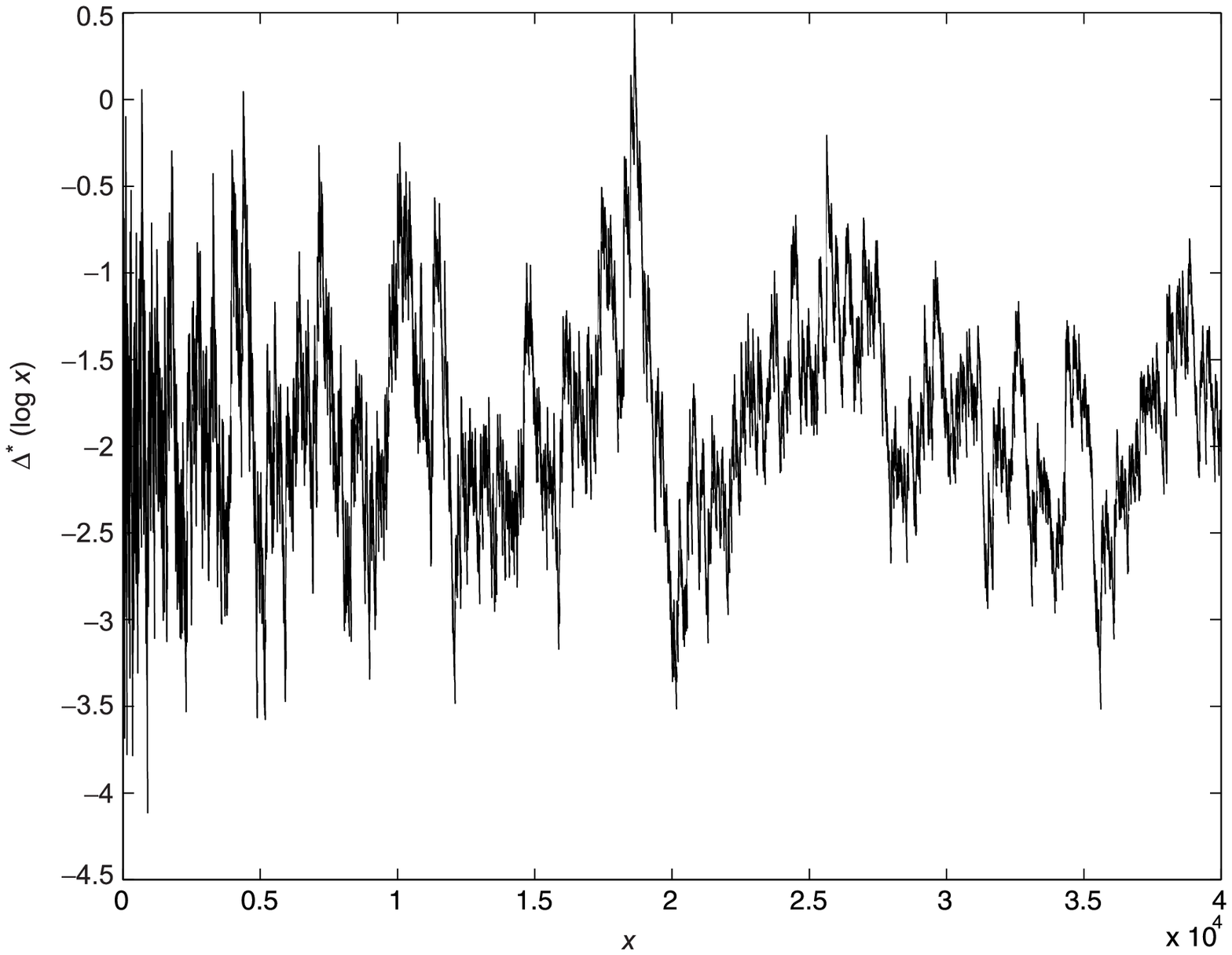}}\vspace{-.8pc}
\caption{}\vspace{.5pc}
\end{figure}

The error is given by the $O$-term in Proposition~4.1 after
approximating $H(R)$ by $3\e^{R}/2$, i.e.,
\begin{equation*}
\hbox{Error}=\hbox{visible}-\frac{3}{2} \e^R+\frac{3}{2} \e^{R/2}.
\end{equation*}

Note that the number of invisible points is relatively small, in
accordance with Proposition~4.1. In fact, following the arguments of its
proof and truncating the series $\sum\mu(n)H(R/n)$ to $n\le 6$, one can
expect\pagebreak
\begin{equation*}
\hbox{Approx.}=\frac{3}{2} (\e^{R/2}+\e^{R/3}+\e^{R/5}-\e^{R/6})
\end{equation*}
to be a good approximation for the number of invisible points.
Table~2 confirms this assertion for the previous data

Finally we show the graph of $\Delta^*(\log x)$ (we are approximating
$H^*(\log x)$ by $H^*(\text{\rm arc\;cosh}(x/2))$, actually) (figure~1). Note
the bias predicted by Corollary~4.2 due to invisible points.

The aspect of this graphic does not differ from the graphics of
normalized error term in classical circle and divisor problem, but in
this case it is not known that there is a limit distribution (cf.
\cite{Hea}).

\section*{Acknowledgements}

I owe my deepest gratitude to E~Valenti, who was fond of non-hyperbolic
orchards.

\end{document}